\documentclass[11pt]{amsart}
\usepackage[french, english]{babel}
\usepackage{amsmath,amssymb,amsthm,mathrsfs, yfonts,url,lscape}
\usepackage{xypic}
\newtheorem{thm}{Theorem}

\newtheorem{prop}[thm]{Proposition}
\newtheorem{cor}[thm]{Corollary}

\newtheorem{rem}[thm]{Remark}

\newtheorem{lem}[thm]{Lemma}

\newtheorem{alg}[thm]{Algorithm}
\newenvironment{prf}{\begin{proof}[\bf{Proof}]}{\end{proof}}

\newcommand{\field}[1]{\mathbb{#1}}
\newcommand{\Q}{\field{Q}}
\newcommand{\C}{\field{C}}
\newcommand{\R}{\field{R}}

\newcommand{\Z}{\field{Z}}
\newcommand{\A}{\field{A}}
\newcommand{\F}{\field{F}}
\newcommand{\G}{\field{G}}

\newcommand{\CO}{\mathcal{O}}

\newcommand{\gp}{\mathfrak{p}}
\newcommand{\gn}{\mathfrak{n}}
\newcommand{\ga}{\mathfrak{a}}
\newcommand{\gb}{\mathfrak{b}}
\newcommand{\gc}{\mathfrak{c}}
\newcommand{\gq}{\mathfrak{q}}
\newcommand{\nr}{\mathrm{nr}}

\DeclareMathOperator{\Gal}{Gal}
\DeclareMathOperator{\GL}{GL}
\DeclareMathOperator{\GSp}{GSp}
\DeclareMathOperator{\GU}{GU}

\title[Algebraic modular forms]{\bf On the computation of algebraic modular\\ forms on compact inner forms of $\mathbf{GSp}_4$}
\author{Lassina Demb\'el\'e}
\address{Warwick Mathematics Institute, University of Warwick, Coventry CV4 7AL, United Kingdom}
\email{l.dembele@warwick.ac.uk}
\date{\today}
\begin{document}
\maketitle

\begin{abstract} In this paper, we describe an algorithm for computing algebraic modular forms on compact inner forms of $\GSp_4$ over totally real number fields. By analogues of the Jacquet-Langlands correspondence for $\GL_2$, this algorithm in fact computes Hecke eigensystems of Hilbert-Siegel modular forms of genus 2. We give some examples of such eigensystems over $\Q(\sqrt{2})$.
\end{abstract}

\section*{Introduction} 
Let $F$ be a totally real number field and $D/F$ a quaternion algebra over $F$ that is ramified at all infinity places of $F$. Let $G/\Q$ be the algebraic group whose $\Q$-rational points are given by the (quaternionic) unitary similitude group $\GU_2(D)$. Then, $G$ is an inner form of $G':=\mathrm{Res}_{F/\Q}(\GSp_4)$ such that $G(\R)$ is compact modulo its center. The Langlands philosohpy predicts that there is a transfer map between automorphic forms (or algebraic modular forms in the sense of Gross~\cite{gross1}) on $G$ and Hilbert-Siegel modular forms on $G'/\Q$. In fact, the existence of such a transfer map was first conjectured by Ihara, Hashimoto-Ibukiyama, and Ibukiyama~\cite{ihara, hashimoto1, ibukiyama} for $F=\Q$, as an analogue in genus 2 of the Jacquet-Langlands correspondence for $\GL_2$ or, more precisely, the Eichler theta correspondence. And recently, Sorensen~\cite{sorensen1} established that such a correspondence does indeed exist when the degree of $F$ is even.

Under the assumption of the existence of a Jacquet-Langlands correspondence between $G$ and $G'$, Cunningham-Demb\'el\'e~\cite{clifton1} give an algorithm for the computation of Hecke eigensystems of Hilbert-Siegel modular forms of genus 2 over $\Q(\sqrt{5})$. Here, we extend the algorithm in that paper to any arbitrary totally real number field $F$. Our approach borrows from Hashimoto-Ibukiyama and Ibukiyama~\cite{hashimoto1, ibukiyama} who provided numerical evidences in support of their conjecture by explicit computations with Brandt matrices. And also from Lansky-Pollack~\cite{lansky} who computed examples of algebraic modular forms on $G$ for $F=\Q$ and determined their Satake parameters. (We note, in passing, that the computations in~\cite{lansky} were used to predict the existence of symmetric cube liftings from $\mathrm{PGL}_2$ to $\mathrm{PGSp}_4$, a result which has now been proved by Ramakrishnan-Shahidi~\cite{ramakrishnan1}.) 

By work of Ghitza~\cite{ghitza1}, all Hecke eigenvalues systems of $\bmod\,p$ Siegel modular forms come from Hecke eigenvalues systems computed from the compact form $G$ obtained from the quaternion algebra $D/\Q$ ramified at $\infty$ and $p\ge 3$ only. His result combined with our algorithm provide a mean to gather numerical evidence on the ``weight part'' of the analogue of the Serre conjecture for $\GSp_4$ (Herzig and Tilouine~\cite{HerzigTilouine}), and we hope to do that. With some changes, this algorithm should also works for other algebraic groups which satisfy the conditions in Gross~\cite[Proposition 1.14]{gross1}, especially compact forms of $\mathrm{U}(3)$. And we also hope to use it for this latter group in order to investigate a conjecture of Bella\"{\i}che-Graftieaux~\cite{bellaiche1} which is an analogue of the so-called Ihara Lemma.

The outline of the paper is as follows. In Section~\ref{sec:innerforms}, we give a brief review of compact inner forms $G$ of $G'$ and the theory of lattices on integral models of $G$. In Section~\ref{sec:algmodforms}, we recall the theory of algebraic modular forms on $G$. Then in Section~\ref{sec:algorithm}, we give some details regarding the implementation of the algorithm. Finally, in Section~\ref{sec:example} we give numerical examples over the real quadratic field $\Q(\sqrt{2})$.

\noindent
{\bf Acknowledgments.} This project started when the author was still a PIMS postdoctorate at the University of Calgary, and completed while he was at the Institut f\"ur Experimentelle Mathematik with an SFB/TR 45 Fellowship of the Deutsche For\-schungsgemeinschaft. He would like to thank these institutions for their financial support. He would also like to thank Clifton Cunningham for his constant support and encouragement. The author would like to thank William Stein for his generosity in allowing him to use the {\sf Sage} cluster in the earlier stages of the project. Finally, the author would like to thank Kevin Buzzard, Steve Donnelly and Dick Gross for helpful email  exchanges.

\section{Inner forms and integral models of $\GSp_4$}\label{sec:innerforms}
Throughout this paper, we let $F$ be a  totally real number field of degree $g$ and  $\mathcal{O}_F$ its ring of integers. For every real embedding $v:\,F\hookrightarrow\R$ and every element $a\in F$, we let $a_v:=v(a)$ be the image of $a$ under $v$. (We recall that an element $a\in F$ is said to be totally positive if $a_v>0$ for all $v\mid\infty$.) Let $D$ be a quaternion algebra over $F$ which is ramified at all infinite places;  we let $\Sigma$ be the set of finite primes at which $D$ is ramified. We denote by $u\mapsto\bar{u}$ and by $\nr: D\to F$ the involution and  the reduced norm of $D$, respectively. We fix a maximal order $\CO_D$ in $D$.  We choose a finite extension $E/F$, contained in $\C$, and such that there is a splitting isomorphism $j:\, D\otimes_\Q E\simeq \mathrm{M}_2(E)^g$.  For any prime $\mathfrak{p}\subset\mathcal{O}_F$, we denote the completions of $F$ and $\mathcal{O}_F$ at $\mathfrak{p}$ by $F_\mathfrak{p}$ and $\mathcal{O}_{F_\mathfrak{p}}$, respectively. Similarly, we denote the completions of $D$ and $\CO_D$ at $\mathfrak{p}$ by $D_{\gp}$ and $\CO_{D_{\gp}}$.  For any finite prime $\mathfrak{p}\notin \Sigma$, we fix an isomorphism $\CO_{D_{\gp}}\simeq\mathrm{M}_2(\mathcal{O}_{F_\mathfrak{p}})$ and extend it to $D_{\mathfrak{p}}\simeq\mathrm{M}_2(F_{\mathfrak{p}})$.  

\subsection{A review on $\CO_D$-lattices in $D^2$}

Here, we recall some results about lattices in $D^2$ and supplement them. We consider $D^2$ with its left $D$-module structure. The endomorphism ring of $D^2$ is $\mathrm{M}_2(D)$, the ring of 2 by 2 matrices with entries in $D$. The involution of $D$ induced an involution on $\mathrm{M}_2(D)$ by action on the entries. Let $A$ be a matrix in $\mathrm{M}_2(D)$. We define the determinant of $A$ as follows. We choose an extension $E/F$ which splits $D$ and let $\det(A)=\det_E(A\otimes 1)$ in $\mathrm{M}_2(D\otimes_FE)\simeq\mathrm{M}_4(E)$. It is not hard to see that $\det(A)$ does not depend on the choice of splitting field. We say that $A$ is hermitian if $A=\bar{A}^t$. In that case, 
$$A=\begin{pmatrix}s&\bar{r}\\ r&t\end{pmatrix}\,\,\mbox{\rm with}\,\,s,t\in F,\,r\in D,$$ and we define $\det_D(A):=st-\nr(r)$. Then, we see that $\det(A)=\det_D(A)^2$.

A quaternion hermitian form on $D^2$ is a map $H:\,D^2\times D^2\to D$ such tat
$$H(u,v)=\overline{H(v,u)}\,\mbox{\rm and}\, H(a u, b v)=aH(u,v)\bar{b}\,\,\forall\,a,b\in D,\,u,v\in D^2.$$
By choosing a basis, we can write $H$ as
$$H(u,v):=u A\bar{v}^t,\,\,\,u,v\in D^2,$$
where $A$ is a hermitian matrix. We say that $A$ is totally positive and write $A>0$ if the form $(u,v)\mapsto uA\bar{v}^t$ is totally positive definite.
 
 From now on, we fix such a totally positive form $H$ on $D^2$ and let $A$ be the associated matrix. A left $\CO_D$-lattice $L\subset D^2$ is an $\mathcal{O}_F$-lattice in $D^2$ which is also a left $\CO_D$-module. As in Shimura~\cite[Section 1]{shimura2}, we define the {\it norm} of such a lattice $L$, with respect to $H$ (or equivalently $A$), as the two-sided $\CO_D$-ideal $\nu_A(L)$ generated by the set $\{H(u,v):\,u,v\in L\}$, and its {\it dual} by
$$L^\#:=\left\{u\in D^2: H(u, L)\subset\CO_D\right\}.$$
 We see that $\nu_A(L)^{-1}L\subseteq L^\#$ since 
$$H(au,v)=aH(u,v)\in\nu_A(L)^{-1}\nu_A(L)=\CO_D\,\,\forall\,a\in\nu_A(L)^{-1},\,u,v\in L.$$ We say that $L$ is  {\it modular} with respect to $H$ (or $A$) if there is an equality, and that $L$ is {\it integral} if the restriction of $H$ to $L\times L$ is $\CO_D$-valued. We recall that a lattice is said to be maximal if its is maximal with respect to the usual inclusion among all lattices of the same norm. We see that if $L$ is integral and maximal then it is self-dual, hence modular, since $\nu_A(L)=\CO_D$. From this and the classification provided by~\cite[Propositions 4.5 and 4.6]{shimura1}, it follows that all the integral maximal lattices belong to the same genus, which we call the {\it principal} genus of $H$. 

Let $L$ and $M$ be two $\CO_D$-lattices in $D^2$. We recall the definition of the {\it index} $[L:M]$ of $M$ in $L$. First, assume that $M\subseteq L$. Then, the quotient $L/M$ is a finite $\CO_D$-module which admits a composition series
$$L/M=L_0\supset L_1\supset \cdots\supset L_r\supset L_{r+1}=0,$$
in which $L_i/L_{i+1}$ is isomorphic to $\CO_D/\mathfrak{m}_i$ for some maximal left $\CO_D$-ideal. By the Jordan-H\"older Theorem, the set ${\mathfrak{m}_i}$ is uniquely determined up to isomorphism by $L/M$. We define the index of $M$ inside $L$ by
$$[L:M]:=\prod_{i=0}^r \mathrm{nr}(\mathfrak{m}_i).$$
In general, one defines the index of $M$ in $L$ as $$[L:M]:=[L: L\cap M][M: M\cap L]^{-1}=[L: P][M:P]^{-1},$$ where we choose any $\CO_D$-lattice $P\subset L\cap M$. For the proof of the following lemma we refer to Coulangeon~\cite[Lemma 2.2.1]{coulangeon}, where one can easily see that the base field $\Q$ plays no special r\^ole.

\begin{lem}\label{lem:index1} Let $L=\mathfrak{a}_1e_1\oplus\mathfrak{a}_2e_2$ be an $\CO_D$-lattice, where $\mathfrak{a}_1$ and $\mathfrak{a}_2$ are frac\-tio\-nal $\CO_D$-ideals. The index of $L$ in $L^\#$ is a square and we have
$$[L^\#:L]=\left(\mathrm{det}_D(H(e_i,e_j))\mathrm{nr}(\mathfrak{a}_1\mathfrak{a}_2)\right)^2.$$
\end{lem}
We define the discriminant $\textfrak{d}_L$ of $L$ to be the square root of $[L^\#:L]$. We then obtain the following lemma. %For the next lemma, we again refer to Coulangeon~\cite[Proposition 2.2.2]{coulangeon}.

\begin{lem}\label{lem:index2} Let $L$ and $M$ be lattices in $D^2$. Then, we have $\textfrak{d}_M=\textfrak{d}_L[L:M]$.
\end{lem}
From Lemmas~\ref{lem:index1} and ~\ref{lem:index2}, it follows immediately that one can equivalently define $[L:M]$ as the $\CO_F$-ideal generated by  the set $\{\mathrm{det}_D(\alpha\bar{\alpha}^t): L\alpha\subseteq M\}$.  And we see that $\textfrak{d}_L=\nu_A(L)^2$ when $L$ is modular.

\begin{lem}\label{lem:index-norm} Let $L$ and $M$ be two $\CO_D$-lattices in $D^2$ which belong to the same genus. Then $\nu_A(M)^2=\nu_A(L)^2[L:M]$.
\end{lem}

\begin{prf} It is enough to prove this equality locally. So, let $\gp\subset\CO_F$ be a prime. By definition of genera, we have $M_\gp=L_\gp g_\gp$ for some $g_\gp\in \GL_2(D_\gp)$ such that $g_\gp A\bar{g}_\gp^t=\nu_A(g_\gp) A$ with $\nu_A(g_\gp)\in F_\gp^\times$. This implies that $\nu_A(M_\gp)=\nu_A(L_\gp)\nu_A(g_\gp)$ and $\mathrm{det}_D(g_\gp\bar{g}_\gp^t)^2=\det(g_\gp\bar{g}_\gp^t)=\nu_A(g_\gp)^4$. Hence, we have
$$[L_\gp:M_\gp]=(\mathrm{det}_D(g_\gp\bar{g}_\gp^t))=(\nu_A(g_\gp))^2=\nu_A(M_\gp)^2\nu_A(L_\gp)^{-2}.$$
\end{prf}

\begin{prop}\label{prop:principal-genus} Let $L$ and $M$ be maximal $\CO_D$-lattices in $D^2$ such that $L$ is modular. Then $M$ belongs to the genus of $L$ if and only if $M$ is modular and the index $[L:M]$ is a square $\CO_F$-ideal.
\end{prop}

\begin{prf} An easy local  calculation shows that modularity is a genus property. But for a modular lattice $M$,  we have $\nu_A(M)^2=\nu_A(L)^2[L:M]$ by Lemmas~\ref{lem:index2} and~\ref{lem:index-norm}. Thus we can write $\nu_A(M)=\nu_A(L)\gb$, where $\gb$ is a two-sided $\CO_D$-ideal. The classification of genera of maximal lattices given by~\cite[Propositions 4.5 and 4.6]{shimura1} then implies that $M$ belongs to the genus of $L$ if and only if $\gb$ is an $\CO_F$-ideal; {\it i.e.}, if and only if $[L:M]$ is a square ideal.
\end{prf}

\begin{cor}\label{cor:principal-genus} Let $L$ be a maximal $\CO_D$-lattice in $D^2$. Then $L$ belongs to the principal genus of $H$ if and only if $L$ is modular and $\textfrak{d}_L$ is a square.
\end{cor}

Let $L$ be a left $\CO_D$-lattice. Since $\CO_D$ is maximal, it is hereditary; hence by Reiner~\cite[Theorem 27.8]{reiner}, there exists a left $\CO_D$-ideal $\mathfrak{a}$ and $\alpha\in\GL_2(D)$ such that $L=(\CO_D\oplus\ga^{-1})\alpha$. Furthermore, by the strong approximation theorem (or~\cite[Theorem 6.14]{shimura2}), we can choose $\ga$ to be two-sided since the class number of (the maximal order) $\mathrm{M}_2(\CO_D)$ is the narrow class number of $F$. The following generalizes Hashimoto-Ibukiyama~\cite[Proposition 22]{hashimoto1}.

\begin{prop}\label{prop:class-reps} Let $L=(\CO_D\oplus \mathfrak{a}^{-1})\alpha$ be a maximal $\CO_D$-lattice in $D^2$, where $\ga$ is a two-sided ideal and $\alpha\in\GL_2(D)$. Then $L$ belongs to the principal genus of $H$ if and only if $\alpha A\bar{\alpha}^t= m B$, where $m\in \nu_A(L)\cap F^{\times+}$ and 
$$B\in\GL_2(D)\cap\begin{pmatrix}\mathcal{O}_F&\bar{\ga}\\ \ga&\nr(\ga)\end{pmatrix}$$
is a totally positive hermitian matrix such that the ideal $\det_D(B)\nr(\ga)^{-1}$ is a square. 
\end{prop}

\begin{prf} We first observe that, by Lemma~\ref{lem:index1}, we have
$$\textfrak{d}_L=\frac{(\mathrm{det}_D(\alpha\bar{\alpha}^t)\mathrm{det}_D(A))}{\nr(\ga)}.$$
Suppose that $L$ belongs to the principal genus of $H$. By Corollary~\ref{cor:principal-genus}, this means that $L$ is modular and $\textfrak{d}_L=\nu_A(L)^2$ is a square ideal or, equivalently, that $L^\#=\nu_A(L)^{-1}L$ and  $\nu_A(L)$ is an $\CO_F$-ideal. This implies that
$$u(a\alpha A\bar{\alpha}^t)\bar{v}^t=(au)(\alpha A\bar{\alpha}^t)\bar{v}^t\in\CO_D\,\,\forall a\in\nu_A(L)^{-1}\cap F,\,u,v\in\CO_D\oplus\ga^{-1}.$$
Therefore, we must have $\alpha A\bar{\alpha}^t=m B$ for some $m\in\nu_A(L)\cap F^{\times +}$ and 
$$B\in\GL_2(D)\cap\begin{pmatrix}\mathcal{O}_F&\bar{\ga}\\ \ga&\nr(\ga)\end{pmatrix}.$$
We see that $B$ is totally positive; and the discriminant formula and the relation $\textfrak{d}_L=\nu_A(L)^2$ show that $\det_D(B)\nr(\ga)^{-1}$ is a square.

Conversely, if there exists a pair $(m, B)$ as above, then $\textfrak{d}_L$ is a square. Thus we only need to show that $L$ is modular. To this end, we write $L=\CO_De_1\oplus\ga^{-1}e_2$ and $L^\#=\CO_De_1^\#\oplus\bar{\ga}e_2^\#$ where the basis $\{e_1, e_2\}$ is given by the rows of $\alpha$. In that basis, the hermitian form is given by the matrix $\alpha A\bar{\alpha}^t=m B$ where $$B=\begin{pmatrix}s&\bar{r}\\ r&t\end{pmatrix}\in\GL_2(D)\cap\begin{pmatrix}\mathcal{O}_F&\bar{\ga}\\ \ga&\nr(\ga)\end{pmatrix}.$$ Using the facts that $m\in\nu_A(L)\cap F^{\times +}$ and $\det_D(B)\in\nr(\ga)$, we see that
\begin{align*}
\left(\alpha A\bar{\alpha}^t\right)^{-1}&=\frac{1}{m\mathrm{det}_D(B)}\begin{pmatrix}t&-\bar{r}\\ -r&s\end{pmatrix}\in\frac{\nu_A(L)^{-1}}{\nr(\ga)}\begin{pmatrix}\nr(\ga)&\bar{\ga}\\ \ga&\CO_F\end{pmatrix}\\
&=\nu_A(L)^{-1}\begin{pmatrix}\CO_F&\ga^{-1}\\ \bar{\ga}^{-1}&\nr(\ga)^{-1}\end{pmatrix}.
\end{align*}
By combining this with the equality $(e_1^\#,e_2^\#)^t=(\alpha A\bar{\alpha}^t)^{-1}(e_1,e_2)^t$, we get the inclusion $L^\#\subset\nu_A(L)^{-1}L$; hence $L$ is modular.
\end{prf}

\begin{rem}\label{rem:class-reps}\rm Proposition~\ref{prop:class-reps} provides a simple criterion to search for lattices which belong to the principal genus of $H$; and in practice we can assume that $m=1$. Indeed, for every $m\in F^{\times +}$,~\cite[Proposition 6.13]{shimura2} shows that there exists $\delta\in\GL_2(D)$ such that $\delta A\bar{\delta}^t=m A$, thus we can replace $\alpha$ by $\alpha\delta^{-1}$ if necessary. Furthermore, it is not hard to see that we can choose $B$ such that $\det_D(B)$ is a $\nr(\ga)$-unit.
\end{rem}

From now on, we will denote the lattice $L=(\CO_D\oplus\ga^{-1})\alpha$ by $L_{\ga,\alpha}$ and its discriminant by $\textfrak{d}_{\ga,\alpha}$.

\subsection{Inner forms of $\GSp_4$ and integral models}\label{subsec:inner-forms}
Let $A$ be a totally positive hermitian matrix. We obtain a reductive group $G^A/\Q$ by putting
$$G^A(R):=\GU_2^A(F\otimes_\Q R)=\left\{\gamma\in\mathrm{M}_2(D\otimes_\Q R)\Big\vert\begin{array}{l}\gamma A\bar{\gamma}^t=\nu_A(\gamma) A\\ \nu_A(\gamma)\in(F\otimes_\Q R)^\times\end{array}\right\},$$
for any $\Q$-algebra $R$. This determines a morphism $\nu_A: G^A\to\G_m$ of al\-ge\-braic groups called the similitude factor. The group $G^A/\Q$ is an inner form of $\mathrm{Res}_{F/\Q}(\GSp_4)$ and, by combining~\cite[Proposition 2.1]{shimura1} and~\cite[Lemma 4]{ghitza1}, we can obtain an explicit isomorphism $G^A(E)\simeq\GSp_4(E)^g$.

Let $L$ be a maximal integral $\CO_D$-lattice in $D^2$. We obtain an integral structure $\underline{G}_L^A$ on $G^A$ over $\Z$ by putting
\begin{align*}
\underline{G}_L^A(R):=&\left\{\gamma\in\mathrm{End}_{\CO_D\otimes_\Z R}(L\otimes_\Z R)\Big\vert\begin{array}{l}\gamma A\bar{\gamma}^t=\nu_A(\gamma)A\\ \nu_A(\gamma)\in(\CO_F\otimes_\Z R)^\times\end{array}\right\}\\
=&\GU^A(L\otimes_\Z R),
\end{align*}
for any $\Z$-algebra $R$. We call the principal genus of $G^A$, the one containing all the integral models obtained from integral maximal lattices in the principal genus of $H$. We observe that, for every prime $p$, 
$$\underline{G}_L^A(\Z_p)=\prod_{\gp\mid p}\GU^A(L\otimes_{\CO_F}\CO_{F_\gp}),$$ where $\GU^A(L\otimes_{\CO_F}\CO_{F_\gp})$ is hyperspecial if and only if $\gp\notin\Sigma$. Therefore, $\underline{G}_L^A$ is an integral model in the sense of Gross~\cite{gross2} if and only if $\Sigma=\emptyset$. This is equivalent to saying that $[F:\Q]$ is even and $D/F$ is the unique quaternion algebra ramified at $v\mid\infty$ only, in which case there is only one genus.

\section{Algebraic modular forms and Hecke action} \label{sec:algmodforms}

From now on, we let $A=\mathbf{1}_2$ and we simply denote the group $G^A$ and its integral model associated to $L=\CO_D^2$ by $G$ and $\underline{G}$, respectively. We recall the theory of algebraic modular forms on $G$ in the sense of Gross~\cite{gross1}. We then introduce the notion of algebraic modular forms on integral models of $G$ and explain how the two relate. 

\subsection{The Hecke module of algebraic modular forms}
For any finite prime $\mathfrak{p}\notin \Sigma$, we choose an isomorphism $\GU_2(\CO_{D_\gp})\cong\GSp_4(\mathcal{O}_{F_\mathfrak{p}})$, which is compatible with the splitting isomorphism $\CO_{D_\gp}=\mathrm{M}_2(\mathcal{O}_{F_\mathfrak{p}})$ we fixed earlier, and we consider the maximal compact open subgroup
$$\underline{G}(\hat{\Z})=\prod_{\mathfrak{p}\in\Sigma}\GU_2(\CO_{D_\gp})\times\prod_{\mathfrak{p}\notin\Sigma}\GSp_4(\mathcal{O}_{F_\mathfrak{p}}).$$
We let $U=\prod_{\mathfrak{p}}U_\mathfrak{p}\subseteq\underline{G}(\hat{\Z})$ be a compact open such that $U_\mathfrak{p}$ is maximal for each prime $\mathfrak{p}\in\Sigma$.

Let $\rho: G\to\GL(V)$ be an irreducible algebraic representation. For any subfield $\Q\subseteq E\subseteq\C$, we denote the realization of $(\rho,\,V)$ on $E$ in the same way. We recall that, when $E$ splits $D$, then this realization is of the form $\rho=\bigotimes_{v\mid\infty}\rho_v$ and $V=\bigotimes_{v\mid\infty}V_v$, where each
$$\rho_v:\,\GSp_4(E)\to\GL(V_v)$$ is an irreducible representation of highest weight $(a_1,a_2; b)\in\Z^3$, with $a_1\ge a_2\ge 3$ and $b\equiv a_1+a_2\mod 2$.

The space of algebraic modular forms on $G(\hat{\Q})$ of weight $V$ and level $U$ is defined by
$$M_{\rho}(U):=\left\{f:\,G(\hat{\Q})/U\to V\big\vert\,f|_{\rho}\gamma=f,\,\,\mbox{\rm for all}\,\,\gamma\in G(\Q)\right\},$$
where $f|_{\rho}\gamma(x)=f(\gamma x)\rho(\gamma)$, for all $x\in G(\hat{\Q})$ and $\gamma\in G(\Q)$.
The Hecke algebra $\mathbf{T}_{\rho}(U)$ acting on this space is the $\Z$-algebra generated by the operators defined as follows. For any $u\in G(\hat{\Q})$, write the finite disjoint union $UuU=\coprod_i u_iU$ and put
$$f|_{\rho}[UuU](x)=\sum_{i}f(xu_i).$$ 
For a prime $\mathfrak{p}\notin\Sigma$, the local Hecke algebra at $\mathfrak{p}$ is generated by the two operators $T_1(\mathfrak{p})$ and $T_2(\mathfrak{p})$ corresponding to the double cosets of the diagonal matrices $\mathrm{diag}(1, 1, \varpi_{\mathfrak{p}},\,\varpi_{\mathfrak{p}})$ and $\mathrm{diag}(1, \varpi_{\mathfrak{p}}, \varpi_{\mathfrak{p}}^2,\,\varpi_{\mathfrak{p}})$ respectively, where $\varpi_{\mathfrak{p}}$ is a uniformizer at $\mathfrak{p}$. Hence the Hecke algebra $\mathbf{T}_{\rho}(U)$ acting on $M_{\rho}(U)$ is the $\Z$-algebra generated by the operators $T_1(\mathfrak{p})$ and $T_2(\mathfrak{p})$ for all $\mathfrak{p}\notin\Sigma$.

When $\rho=\mathrm{triv}$ is the trivial representation, we let $I_{\rho}(U)$ be the subspace of $M_{\rho}(U)$ given by
$$I_{\rho}(U):=\left\{f:\,G(\Q)\backslash G(\hat{\Q})/U\stackrel{\nu}{\longrightarrow} F^{\times+}\backslash (\A_F^f)^\times/\nu(U)\longrightarrow \C \right\}.$$
Then, we define
\begin{eqnarray*}
S_{\rho}(U):=\left\{\begin{array}{ll}
M_{\rho}(U),&\rho\neq\mathrm{triv}\\
&\\
M_{\rho}(U)/I_{\rho}(U),&\rho=\mathrm{triv}.
\end{array}\right.
\end{eqnarray*} The action of Hecke preserves $I_{\rho}(U)$, and thus induces an action of $\mathbf{T}_{\rho}(U)$ on $S_{\rho}(U)$. We refer to~Ihara, Hashimoto-Ibukiyama, and Ibukiyama~\cite{ihara, hashimoto1, ibukiyama} for the  conjectural analogue of the Jacquet-Langlands correspondence for $\GSp_4$ (over $F=\Q$) which relates the Hecke module $S_\rho(U)$ to spaces of Hilbert-Siegel modular forms, and in particular to Sorensen~\cite{sorensen1} who shows that this conjecture is true when the degree of $F$ is even.

\subsection{Integral models of $G$ and the Hecke module $M_{\rho}(U)$}
 
Let $\mathrm{Cl}(\underline{G})$ be a complete set of representatives of the conjugacy classes in the genus of $\underline{G}$ as described in~\cite{gross2}. There is a canonical bijection between $\mathrm{Cl}(\underline{G})$ and the double coset $G(\Q)\backslash G(\hat{\Q})/\underline{G}(\hat{\Z})$, which also parametrizes isomorphism classes of left $\CO_D$-lattices in the genus of $\CO_D^2$. By abuse of notation, we let $(\ga,\alpha)$ be the integral model attached to the lattice $L_{\ga,\alpha}$. Let $\Sigma_U:=(\Sigma_U^1,\Sigma_U^2)$ be a pair of sets consisting of primes such that $\Sigma_U^1\cap\Sigma_U^2=\emptyset$ and each $\Sigma_U^i$ is disjoint from $\Sigma$ and the set of primes at which $U_\mathfrak{p}$ is not maximal; further assume that $\Sigma_U^1$ generates the narrow class group $\mathrm{Cl}^+(F)$. By the strong approximation theorem we can choose the representatives  $(\ga,\alpha)\in \mathrm{Cl}(\underline{G})$ to be supported in $\Sigma_U$; {\it i.e.}, such that $\mathrm{nr}(\mathfrak{a})$ is supported in $\Sigma_U^1$ and $\alpha$ is $\Sigma_U^2$-integral. (In practice, one can choose the set $\Sigma_U^2$ to be a singleton, but it can't be empty unless the class number of $\underline{G}$ is equal to $h^+:=\#\mathrm{Cl}^+(F)$.)

\medskip
We define $\mathcal{F}_{\underline{G}}(U):=\underline{G}(\hat{\Z})/U$ (this is  a finite product of finite sets).  Let $(\mathfrak{a},\alpha)\in\mathrm{Cl}(\underline{G})$ be a representative; the group of automorphisms of the corresponding lattice $L_{\mathfrak{a},\alpha}$ is given by
\begin{align*}
\Gamma_{\mathfrak{a},\alpha}:=\mathrm{Stab}_{G(\Q)}(L_{\mathfrak{a},\alpha})=\left\{\gamma\in G(\Q)\cap \alpha^{-1}\begin{pmatrix}\CO_D&\mathfrak{a}^{-1}\\ \mathfrak{a}&\CO_D\end{pmatrix}\alpha:\,\nu(\gamma)\in\mathcal{O}_F^{\times +}\right\},
\end{align*}
and the quotient $\overline{\Gamma}_{\ga,\alpha}=\Gamma_{\ga,\alpha}/\CO_F^\times$ is finite. We define
$$M_{\rho}^{\Gamma_{\mathfrak{a},\alpha}}(U):=\left\{f:\,\mathcal{F}_{\underline{G}}(U)\to V:\,\,f|_{\rho}\gamma=f\,\,\forall\,\gamma\in\Gamma_{\mathfrak{a},\alpha}\right\}.$$
Let $(\mathfrak{b},\beta)\in\mathrm{Cl}(\underline{G})$ be another representative, and recall that
\[
\mathrm{Hom}_{\CO_D}(L_{\mathfrak{a},\alpha},\,L_{\mathfrak{b},\beta}):=\left\{\gamma\in\mathrm{M}_2(D)\Big|\,L_{\mathfrak{a},\alpha}\gamma\subseteq L_{\mathfrak{b},\beta}\right\}=\alpha^{-1}\begin{pmatrix}\CO_D&\mathfrak{b}^{-1}\\ \mathfrak{a}& \mathfrak{a}\mathfrak{b}^{-1} \end{pmatrix}\beta.
\]
Now, let $\mathfrak{p}\notin \Sigma$ be a prime ideal in $\mathcal{O}_F$, and choose $a'\in\mathfrak{a}^{-1}$ and $b\in\mathfrak{b}$ such that $\mathfrak{a}_\mathfrak{p}^{-1}=(a_\mathfrak{p}')$ and $\mathfrak{b}_\mathfrak{p}=(b_\mathfrak{p})$; (this is equivalent to saying that
$v_{\mathfrak{p}}(\mathrm{nr}(a'))=-v_{\mathfrak{p}}(\mathrm{nr}(\mathfrak{a}))$ and $v_{\mathfrak{p}}(\mathrm{nr}(b))=v_{\mathfrak{p}}(\mathrm{nr}(\mathfrak{b}))$).
Then, for every element $\gamma\in\mathrm{Hom}_{\CO_D}(L_{\mathfrak{a},\alpha},\,L_{\mathfrak{b},\beta})$, the matrix
$\mathrm{diag}(1,\,a')\alpha\gamma\beta^{-1}\mathrm{diag}(1,\,b)$ is $\mathfrak{p}$-integral. Thus, we get a well-defined map
\begin{eqnarray*}
\vartheta_{(\mathfrak{a},\alpha),(\mathfrak{b},\beta)}:\,\mathrm{Hom}_{\CO_D}(L_{\mathfrak{a},\alpha},\,L_{\mathfrak{b},\beta})&\to&\mathrm{M}_2(\CO_{D_\gp})\\
\gamma&\mapsto&\begin{pmatrix}1& 0\\ 0&a_\mathfrak{p}'\end{pmatrix}\alpha_\mathfrak{p}\gamma_\mathfrak{p}\beta_\mathfrak{p}^{-1}\begin{pmatrix}1&0\\ 0& b_\mathfrak{p}\end{pmatrix}
\end{eqnarray*}
Let $\F_\mathfrak{p}=\mathcal{O}_F/\mathfrak{p}$ be the residue field at $\mathfrak{p}$ and $\CO_{D_\gp}\simeq\mathrm{M}_2(\mathcal{O}_{F_\mathfrak{p}})\to\mathrm{M}_2(\F_\mathfrak{p})$ a reduction map (obtained from a fixed splitting map). Then, by composition with $\vartheta_{(\mathfrak{a},\alpha),(\mathfrak{b},\beta)}$, we get a well-defined map
\begin{eqnarray*}
\mathrm{Hom}_{\CO_D}(L_{\mathfrak{a},\alpha},\,L_{\mathfrak{b},\beta})&\to&\mathrm{M}_4(\F_\mathfrak{p})\\
\gamma&\mapsto&\tilde{\gamma}
\end{eqnarray*} 
Next, for $i=1$ and $2$, we define the sets
\begin{eqnarray*}
\lefteqn{\Theta^{\Sigma_U}_i(\mathfrak{p};(\mathfrak{a},\alpha),(\mathfrak{b},\beta)):=}\\
&&\qquad\Gamma_{\mathfrak{a},\alpha}\backslash\left\{\gamma\in G(\Q)\cap\mathrm{Hom}_{\CO_D}(L_{\mathfrak{a},\alpha},\,L_{\mathfrak{b},\beta})\Big\vert\begin{array}{l} [L_{\mathfrak{b},\beta}:L_{\mathfrak{a},\alpha}\gamma]=\mathfrak{p}^i,\\ \mathrm{rank}_{\F_\mathfrak{p}}(\tilde{\gamma})=3-i\end{array}\right\}
\end{eqnarray*}
where we let $\Gamma_{\mathfrak{a},\alpha}$ act on the left. By unravelling definitions, we get 
\begin{eqnarray*}
\lefteqn{\Theta^{\Sigma_U}_1(\mathfrak{p};(\mathfrak{a},\alpha),(\mathfrak{b},\beta)):=}\\
&&\qquad\Gamma_{\mathfrak{a},\alpha}\backslash\left\{\gamma\in G(\Q)\cap\alpha^{-1}\begin{pmatrix}\CO_D&\mathfrak{b}^{-1}\\ \mathfrak{a}& \mathfrak{a}\mathfrak{b}^{-1} \end{pmatrix}\beta\Big\vert \,\nu(\gamma)\textfrak{d}_{\mathfrak{a},\alpha}\textfrak{d}_{\mathfrak{b},\beta}^{-1}=\mathfrak{p}\right\}
\end{eqnarray*}
\begin{eqnarray*}
\lefteqn{\Theta^{\Sigma_U}_2(\mathfrak{p};(\mathfrak{a},\alpha),(\mathfrak{b},\beta)):=}\\
&&\qquad\Gamma_{\mathfrak{a},\alpha}\backslash\left\{\gamma\in G(\Q)\cap\alpha^{-1}\begin{pmatrix}\CO_D&\mathfrak{a}^{-1}\\ \mathfrak{a}& \mathfrak{a}\mathfrak{b}^{-1} \end{pmatrix}\beta\Big\vert\begin{array}{l}\nu(\gamma)\textfrak{d}_{\mathfrak{a},\alpha}\textfrak{d}_{\mathfrak{b},\beta}^{-1}=\mathfrak{p}^2,\\ \mathrm{rank}_{\F_\mathfrak{p}}(\tilde{\gamma})=1\end{array}\right\}.
\end{eqnarray*}
(It is not hard to see that these sets do not depend on the choice of $a'\in\mathfrak{a}^{-1}$, $b\in\mathfrak{b}$ and the splitting map at $\mathfrak{p}$.) We define the linear map
\begin{eqnarray*}
T_i(\mathfrak{p}; (\mathfrak{a},\alpha),(\mathfrak{b},\beta)):\,\,
M_{\rho}^{\Gamma_{\mathfrak{b},\beta}}(U) &\to& M_{\rho}^{\Gamma_{\mathfrak{a}, \alpha}}(U)\\
f&\mapsto&\sum_{u\in\Theta_i^{\Sigma_U}(\mathfrak{p};\,(\mathfrak{a},\alpha),(\mathfrak{b},\beta))}f|_{\rho}u,\,\,i=1,2.
\end{eqnarray*}
The following proposition is proved  in the same way as~\cite[Theorem 2]{dembele2}.
\begin{prop}\label{prop1}
There is an isomorphism of Hecke modules
\begin{eqnarray*}
M_{\rho}(U)\stackrel{\sim}{\to} \bigoplus_{(\mathfrak{a},\alpha)\in \mathrm{Cl}(\underline{G})} M_{\rho}^{\Gamma_{\mathfrak{a},\alpha}}(U),
\end{eqnarray*} 
where the action of the Hecke operator $T_i(\mathfrak{p})$ on the right is given by the collection of linear maps $T_i(\mathfrak{p}; (\mathfrak{a},\alpha),(\mathfrak{b},\beta)),\,(\mathfrak{a},\alpha),(\mathfrak{b},\beta)\in \mathrm{Cl}(\underline{G})$.
\end{prop}
Proposition~\ref{prop1} also gives a decomposition $\displaystyle{S_{\rho}(U)\simeq \bigoplus_{(\mathfrak{a},\alpha)\in \mathrm{Cl}(\underline{G})} S_{\rho}^{\Gamma_{\mathfrak{a},\alpha}}(U)}$.

\section{The algorithm}\label{sec:algorithm}
We now explain how to implement the discussion of Section~\ref{sec:algmodforms} into a concrete algorithm. As with our presentation, the algorithm we present here parallels the one in~\cite{clifton1} and~\cite{dembele2}.

\subsection{Representing hermitian quadratic forms}\label{subsec:forms-rep}
In this subsection, we discuss the problem of representing hermitian quadratic forms by hermitian quadratic forms. This problem arises naturally when testing for equivalence in $\mathrm{Cl}(\underline{G})$ or during the computation of the Hecke action.  We remark that analogous questions appear when dealing with Fourier expansions of Siegel modular forms. 

Let $A$ be a totally positive hermitian matrix and $Q$ the associated quadratic form. The matrix $A$ also gives a vector-valued $D$-hermitian form
\begin{eqnarray*}
\mathrm{M}_2(D)\times\mathrm{M}_2(D)&\to&\mathrm{M}_2(D)\\
(\gamma_1,\gamma_2)&\mapsto&\gamma_1 A\bar{\gamma}_2^t.
\end{eqnarray*} Let $\mathfrak{a}_1$, $\mathfrak{b}_1$, $\mathfrak{c}_1$ and $\mathfrak{d}_1$ be left $\CO_D$-ideals in $D$ so that $\Lambda:=\begin{pmatrix}\mathfrak{a}_1&\mathfrak{b}_1\\ \mathfrak{c}_1&\mathfrak{d}_1\end{pmatrix}$ is an $\CO_D$-lattice in $\mathrm{M}_2(D)$.  Given a hermitian matrix $\eta\in\mathrm{M}_2(D)$, we want to compute the set
\begin{eqnarray*}
\Theta_A(\eta):=\left\{\gamma\in\Lambda:\,\gamma A\bar{\gamma}^t=\eta=\begin{pmatrix}s&\bar{r}\\ r&t\end{pmatrix}\right\}.
\end{eqnarray*}
Let $\gamma=\begin{pmatrix}u\\ v\end{pmatrix}\in\mathrm{M}_2(D)$ be given in rows. Then, it is easy to see that
 $$\gamma A\bar{\gamma}^t=\begin{pmatrix}s&\bar{r}\\ r&t \end{pmatrix}\,\,\Longleftrightarrow\,\,\left\{\begin{array}{lr} Q(u)=s&(1)\\ \\
v A \bar{u}^t=r& (2)\\ \\
Q(v)=t&(3)\end{array}\right.$$
Thus computing $\Theta_A(\eta)$ is the same as representing $s,t\in F$ by $Q$ under the constraint $(2)$. This implies that $\Theta_A(\eta)$ is finite since $A>0$.  By Shimura~\cite[Proposition 2.1]{shimura1}, there is a diagonal matrix $\eta'\in\mathrm{M}_2(F)$ and $P\in\GL_2(D)$ such that
$\eta=P\eta'\bar{P}^t$. By making the change of variables $\gamma=P\gamma'$ and letting $\Lambda=P\Lambda'$, we see that this is again the same as computing the set
\begin{eqnarray*}
\Theta_A(\eta'):=\left\{\gamma'\in\Lambda':\,\gamma' A\bar{\gamma}'{}^t=\eta'=\begin{pmatrix}s'&0\\ 0&t' \end{pmatrix}\right\},
\end{eqnarray*} where the constraint $(2)$ is now an orthogonality condition. Let $\Lambda'[1]$ and $\Lambda'[2]$ be the rows of $\Lambda'$, which are $\CO_D$-lattices in $D^2$. To compute $\Theta_A(\eta')$, we first determined which of the row lattices $(\Lambda'[1], Q)$ or $(\Lambda'[2],Q)$ is better to work with by comparing their determinants. Without loss of generality, assume this would be $\Lambda'[1]$. We then apply the algorithm in Demb\'el\'e-Donnelly~\cite[Subsection 2.2]{dd1} to efficiently compute the set of elements in $(\Lambda'[1],Q)$ which represent $s'$. 

Now, let $u=(a,b)\in\Lambda'[1]$ be a solution to $Q(u)=s'$. We replace $u$ into $(2)$ to get a homogenous linear system in terms of the components of $v=(c,d)\in\Lambda'[2]$.  Solving this linear system in terms of $c$ or $d$ and then replacing into $(3)$, we must now represent $t'$ by a quadratic form $Q_u$ in one variable over an $\CO_D$-ideal. From each solution of $Q_u(w)=t'$, we get $v\in\Lambda'[2]$ which, combined with $u$, gives a matrix $\gamma'\in\Theta_A(\eta')$. (The trick of using the orthogonality condition $(2)$ to cut down the number of variables in $(3)$ is substantially more efficient than the na\"{\i}ve approach of solving $(1)$ and $(3)$ separately and then checking $(2)$ for compatibility.)

\subsection{The representatives $(\mathfrak{a},\alpha)\in \mathrm{Cl}(\underline{G})$}

We explain how to compute representatives for the equivalence classes in the genus of $\underline{G}$.  In order to do so we need the following two lemmas both adapted from~\cite{hashimoto1}.

\begin{lem}\label{lem:equivalence} Let  $L_{\mathfrak{a},\alpha}$ and $L_{\mathfrak{a},\alpha'}$ be in the genus of $\underline{G}$. Then, $L_{\ga,\alpha}$ and $L_{\ga,\alpha'}$ are equivalent if and only if there exist $n\in F^{\times +}$ and $\beta\in\GL(\CO_D\oplus\ga^{-1})$ such that $\alpha'\bar{\alpha}'{}^t=n\beta(\alpha\bar{\alpha}^t)\bar{\beta}^t$.
\end{lem}

\begin{prf} This easily follows from Proposition~\ref{prop:class-reps}.\end{prf}

Let $\ga$ be a two-sided ideal. Then, motivated by Lemma~\ref{lem:equivalence}, we say that two hermitian matrices $\gamma,\,\gamma'\in\GL_2(D)\cap\begin{pmatrix}\CO_F&\bar{\ga}\\ \ga&\nr(\ga)\end{pmatrix}$ are $\ga$-{\it equivalent} if there exist $n\in F^{\times +}$ and $\beta\in\GL(\CO_D\oplus\ga^{-1})$ such that $\gamma'=n\beta\gamma\bar{\beta}^t$.

\begin{rem}\rm We test whether two lattices $L_{\ga,\alpha}$ and $L_{\gb,\beta}$ in the genus of $\underline{G}$ are equivalent by first reducing to the situation in Lemma~\ref{lem:equivalence}. By the observation preceding Proposition~\ref{prop:class-reps},  this amounts to testing whether $\ga$ and $\gb$ are in the same narrow class. In fact, Proposition~\ref{prop:class-reps} and Lemma~\ref{lem:equivalence} imply that the equivalence class of the lattice $L_{\mathfrak{a},\alpha}$ is completely determined by the narrow class of $\mathfrak{a}$, and the $\ga$-equivalence class of $\gamma:=\alpha\bar{\alpha}^t\bmod F^{\times +}$. (See Remark~\ref{rem:class-reps}, and~\cite[Lemma 4.4]{shimura1} by which we can always find a matrix $\alpha\in\GL_2(D)$ such that $\alpha\bar{\alpha}^t=\gamma$.)
\end{rem}

\begin{lem} Let $\mathfrak{a}$ be a two-sided $\CO_D$-ideal and $s\in\mathcal{O}_F^+$, and let 
$$\gamma=\begin{pmatrix}s&\bar{r}\\ r&t\end{pmatrix}\in\GL_2(D)\cap\begin{pmatrix}\CO_F&\bar{\ga}\\ \ga&\nr(\ga)\end{pmatrix}.$$
Then, the $\ga$-equivalence class of $\gamma$  depends on the class $r\bmod {s\ga}$ only.
\end{lem}
To compute the set $\mathrm{Cl}(\underline{G})$, we use the following algorithm.

\medskip
\begin{alg}[]\label{alg:classes-rep}\rm Given a finite set of primes $\Sigma$ and an integral ideal $\gn$ supported outside of $\Sigma$, this returns a set $\mathrm{Cl}(\underline{G})$ of representatives $(\ga,\alpha)$ supported outside of $\Sigma$ and $\gn$.

\begin{enumerate}
\item Initialize $R_{\underline{G}}=\emptyset$ and $M_{\underline{G}}=0$, and compute the mass $\mathrm{Mass}(\underline{G})$.

\item Find a set $\Sigma'$ of prime ideals that generate $\mathrm{Cl}^+(F)$ and such that $\Sigma'$ is disjoint from $\Sigma\cup\{\gp:\gp\mid\gn\}$. Choose a prime $\gq\notin\Sigma'\cup\Sigma\cup\{\gp:\gp\mid\gn\}$.

\item Compute a set $S_{D}$ of integral  two-sided $\CO_D$-ideals supported in $\Sigma'$ and such that, for each $\gb\in\mathrm{Cl}^+(F)$, there is a unique $\ga\in S_D$ such that $[\nr(\ga)]=\gb^2$.  

\item For each $\ga\in S_D$, find a set of representatives $U_\ga$ for the quotient of
$$\left\{u\in\nr(\ga):\, u\,\mbox{\rm is a}\,\nr(\ga)\mbox{-unit and}\,\frac{(u)}{\nr(\ga)}\,\mbox{\rm is a square ideal}\right\}$$ modulo squares. (This is equivalent to a Selmer group computation.)

\item Choose $s\in \mathcal{O}_{F}^+$ (ordered by norm).

\item Choose $\mathfrak{a}\in S_D$ and compute the set  $\mathfrak{a}/s\mathfrak{a}$. 

\item For each $(u,r)\in U_\mathfrak{a}\times(\mathfrak{a}/s\mathfrak{a})$ such that $\mathrm{nr}(r)+u\in \nr(\mathfrak{a})s$,  put $t:=(\mathrm{nr}(r)+ u)/s$. Check if the $\ga$-equivalence class of $\gamma:=\begin{pmatrix}s&\bar{r}\\ r&t\end{pmatrix}$ is already represented in $R_{\underline{G}}$. If not, solve the equation $\alpha\bar{\alpha}^t=\gamma$ for a matrix $\alpha$, which is $\gq$-integral, and compute the finite group $\overline{\Gamma}_{\ga,\alpha}$. Append $(\ga,\alpha)$ to $R_{\underline{G}}$, set $M_{\underline{G}}=M_{\underline{G}}+\frac{1}{\#\overline{\Gamma}_{\mathfrak{a},\alpha}}$. 

\item If  $M_{\underline{G}}=\mathrm{Mass}(\underline{G})$, then $\mathrm{Cl}(\underline{G})=R_{\underline{G}}$, hence stop. Else, go to Step 6. If Step 6 fails, go to Step 5.
\end{enumerate}
\end{alg}

\begin{rem}\rm We test the equivalence of two $\CO_D$-lattices $L$ and $L'$ by using Lemma~\ref{lem:equivalence} or by comparing their theta series. The latter is often enough to quickly distinguish most isomorphism classes and we only use Lemma~\ref{lem:equivalence} for the remaining cases. 

Let $L$ be a modular lattice so that $\nu(L)$ is an $\CO_F$-ideal, and let $\mathfrak{c}$ be a representative of the narrow class of $\nu(L)^{-1}$. We choose $\mu\in F^{\times +}$ such that $\nu(L)\mathfrak{c}=(\mu)$. For every element $(x,y)\in L$, we define its scaled norm with respect to $L$ and $\mathfrak{c}$ by
$$\nu_{L,\mathfrak{c}}(x,y):=\frac{\nu(x,y)}{\mu}.$$
For a given representative $\mathfrak{c}$, this is well-defined up to a totally positive unit. For every integral ideal $\mathfrak{m}\subseteq\mathcal{O}_F$, we define
\begin{align*}
R_{L}(\mathfrak{m}):=\Gamma_L\backslash\left\{(x,y)\in L: \frac{\nu(x,y)}{\nu(L)}=\mathfrak{m}\right\},
\end{align*} 
where $\Gamma_L:=\mathrm{Aut}_{G(\Q)}(L)$. Then $R_L(\mathfrak{m})$ is non-empty if and only if $\nu(L)^{-1}$ and $\mathfrak{m}$ are in the same narrow class. In that case, let $u$ be a totally positive generator of  $\nu(L)\mathfrak{m}$.  
Then the set of elements $(x,y)\in L$ such that $\nu(x,y)=u$ is clearly finite. Thus the quotient $R_{L}(\mathfrak{m})$ is finite. Let $r_{L}(\mathfrak{m}):=\#R_{L}(\mathfrak{m})$ be its cardinality. For each $\tau\in\mathfrak{H}^g$, where $\mathfrak{H}$ is the Poincar\'e upper-half plane, we define
\begin{align*}
\Theta_{L,\mathfrak{c}}(\tau):=&\sum_{(x,y)\in L/\Gamma_L} e^{2\pi i \mathrm{Tr}\left(\nu_{L,\mathfrak{c}}(x,y)\tau\right)}
=\sum_{\mathfrak{m}\subseteq\mathcal{O}_F}^{\qquad\prime} r_{L}(\mathfrak{m}) e^{2\pi i\mathrm{Tr}\left(\nu_{\mathfrak{c}}(\mathfrak{m})\tau\right)},
\end{align*} where the second summation is restricted to all $\mathfrak{m}\subseteq\CO_F$ that belong to the narrow class of $\mathfrak{c}$, and for each such $\mathfrak{m}$ we let $\nu_{\mathfrak{c}}(\mathfrak{m})$ be a totally positive generator of $\mathfrak{m}\mathfrak{c}^{-1}$. A similar argument as in Eichler~\cite[\& 4, Theorem 1]{eichler1} yields that $\Theta_{L,\gc}$ is a Hilbert modular theta series of parallel weight 2 and whose level depends on $\mathfrak{c}$ and the primes in $\Sigma$. 
\end{rem}
 
\subsection{The stabilizers $\Gamma_{\mathfrak{a},\alpha}$ and the sets $\Theta^{\Sigma_U}_i(\mathfrak{p};(\mathfrak{a},\alpha),(\mathfrak{b},\beta))$}\label{subsec:stabilizers}
In this subsection, we describe how the stabilizers and the prelimimary data required for the Hecke action are computed.  To this end, let $(\mathfrak{a},\alpha),(\mathfrak{b},\beta)\in\mathrm{Cl}(\underline{G})$ and $u\in F$ totally positive. We can  compute the finite set
\begin{align*}
S(u; (\mathfrak{a},\alpha),(\mathfrak{b},\beta)):=&\left\{\gamma\in G(\Q)\big\vert\,L_{\mathfrak{a},\alpha}\gamma\subseteq L_{\mathfrak{b},\beta}\,\mbox{\rm and}\,\nu(\gamma)=u\right\}\\
=&\left\{\gamma\in G(\Q)\cap\alpha^{-1}\begin{pmatrix}\CO_D&\mathfrak{b}^{-1}\\ \mathfrak{a}& \mathfrak{a}\mathfrak{b}^{-1} \end{pmatrix}\beta\big\vert\,\gamma\bar{\gamma}^t=u\mathbf{1}_2\right\},
\end{align*}
by making use of the algorithm in Subsection~\ref{subsec:forms-rep} where we first write out:
$$\alpha^{-1}\begin{pmatrix}\CO_D&\mathfrak{b}^{-1}\\ \mathfrak{a}& \mathfrak{a}\mathfrak{b}^{-1} \end{pmatrix}\beta=\begin{pmatrix}\mathfrak{a}_1&\mathfrak{b}_1\\ \mathfrak{c}_1&\mathfrak{d}_1\end{pmatrix}.$$

To obtain $\Gamma_{\mathfrak{a},\alpha}$ (or equivalently $\overline{\Gamma}_{\mathfrak{a},\alpha}$), we simply set $(\mathfrak{a},\alpha)=(\mathfrak{b},\beta)$ and let $U_F^+$ be a set of representatives for $\mathcal{O}_{F}^{\times +}/\left(\mathcal{O}_{F}^{\times}\right)^2$. Then, we get
$$\overline{\Gamma}_{\mathfrak{a},\alpha}=\coprod_{u\in U_F^+}S(u; (\mathfrak{a},\alpha),(\mathfrak{a},\alpha)).$$
For the computation of the representatives of the Hecke double cosets, we first observe that, for $i=1,2$, the set $\Theta^{\Sigma_U}_i(\mathfrak{p};(\mathfrak{a},\alpha),(\mathfrak{b},\beta))$ is non-empty if and only if the ideals $\textfrak{d}_{\mathfrak{a},\alpha}\textfrak{d}_{\mathfrak{b},\beta}^{-1}$ and $\mathfrak{p}^i$ are in the same narrow class. Letting $u$ be a totally positive generator of $\textfrak{d}_{\mathfrak{b},\beta}\textfrak{d}_{\mathfrak{a},\alpha}^{-1}\mathfrak{p}^i$, we then get
\begin{align*}
\Theta^{\Sigma_U}_i(\mathfrak{p};(\mathfrak{a},\alpha),(\mathfrak{b},\beta))&=\Gamma_{\mathfrak{a},\alpha}\backslash S(u; (\mathfrak{a},\alpha),(\mathfrak{b},\beta)). 
\end{align*}

\begin{rem}\rm For each prime $\gp\subset\CO_F$, the cardinalities of these sets satisfy the following identity (see~\cite[Remark 3]{dd1}): for $i=1,2$ and $(\mathfrak{b},\beta)\in\mathrm{Cl}(\underline{G})$,
\begin{align*}
\sum_{(\mathfrak{a},\alpha) \in \mathrm{Cl}(\underline{G})} \#\Theta^{\Sigma_U}_i(\mathfrak{p};(\mathfrak{a},\alpha),(\mathfrak{b},\beta))&=\mathrm{N}\mathfrak{p}^{i-1}(\mathrm{N}\mathfrak{p}+1)(\mathrm{N}\mathfrak{p}^2+1)=\deg T_i(\gp).
\end{align*}
\end{rem}

\subsection{The description of the flag spaces $\mathcal{F}_{\underline{G}}(U)$}

In this subsection, we describe the flag $\mathcal{F}_{\underline{G}}(U)$. The level structure $U$ will be one of the following types: Siegel, Klingen or Borel parahoric (at primes $\mathfrak{p}\notin\Sigma$).  Let $\mathfrak{n}\subseteq\mathcal{O}_F$ be an ideal such that $\mathfrak{p}\nmid\mathfrak{n}$ for $\mathfrak{p}\in\Sigma$. 

\medskip
$\bullet$ {\bf The Siegel parahoric}\\
$$U:=\left\{\gamma\in \underline{G}(\hat{\Z}):\,\gamma\equiv\begin{pmatrix}\ast&\ast&\ast&\ast\\ \ast&\ast&\ast&\ast\\ 0&0&\ast& \ast\\ 0&0&\ast&\ast\end{pmatrix}\mod\mathfrak{n}\right\}.$$ 
We recall the following description from~\cite{clifton1}. Consider the free rank $4$ $\left(\mathcal{O}_{F}/\mathfrak{n}\right)$-module $E:=\left(\mathcal{O}_{F}/\mathfrak{n}\right)^4$ endowed with the symplectic pairing $\langle\,,\,\rangle$
given by the matrix
$$J_2=\begin{pmatrix}0&\mathbf{1}_2\\ -\mathbf{1}_2&0 \end{pmatrix},$$ where $\mathbf{1}_2$ is the identity matrix in
$\mathrm{M}_2(\mathcal{O}_{F}/\mathfrak{n})$. We fix the canonical symplectic basis $e_1,e_2,f_1,f_2$ such that $\langle e_i,\,e_j\rangle=\langle f_i,\,f_j\rangle=0$ and $\langle e_i,\,f_j\rangle=\delta_{ij}$ (Kronecker delta).  Let $M$ be a rank 2 $\left(\mathcal{O}_{F}/\mathfrak{n}\right)$-submodule of $E$. We say that $M$ is {isotropic} if $\langle u,\,v\rangle =0$ for all $u,\,v\in M$. Let  $\mathcal{F}_{\underline{G}}(\mathfrak{n})$ be the set of all rank 2 $\left(\mathcal{O}_{F}/\mathfrak{n}\right)$-submodules of $E$ which are isotropic direct factors. We recall that $\underline{G}(\hat{\Z})$ acts transitively  on $\mathcal{F}_{\underline{G}}(\mathfrak{n})$, and that the stabilizer of the submodule generated by $e_1$ and $e_2$ is the Siegel parahoric $U$. In other words, $\mathcal{F}_{\underline{G}}(\mathfrak{n})$ is the set of $\left(\mathcal{O}_{F}/\mathfrak{n}\right)$-rational points of the Lagrange scheme. 

By using Pl\"ucker's coordinates, we can identify $\mathcal{F}_{\underline{G}}(\mathfrak{n})$ with the set of  $\left(\mathcal{O}_{F}/\mathfrak{n}\right)$-rational points of the Klein quadric as a closed subspace of $\mathbf{P}^5(\mathcal{O}_{F}/\mathfrak{n})$. We then represent each element in $\mathcal{F}_{\underline{G}}(\mathfrak{n})$ by choosing a point $x=(a_0:\cdots:a_5)=[u\wedge v]\in\mathbf{P}^5(\mathcal{O}_{F}/\mathfrak{n})$ such that the submodule $M$ generated by $u$ and $v$ is a rank 2 isotropic direct factor of $E$. We observe that the module $M$ is a direct factor of $E$ if and only if  the ideal $(a_0,\ldots,a_5)$ generated by its coordinates is the unit ideal.

The cardinality of $\mathcal{F}_{\underline{G}}(\mathfrak{n})$ is given by~\cite[Lemma 3.2]{clifton1}:
$$\#\mathcal{F}_{\underline{G}}(\mathfrak{n})=\mathrm{N}(\mathfrak{n})^3\prod_{\mathfrak{p}\mid\mathfrak{n}}\left(1+\frac{1}{\mathrm{N}(\mathfrak{p})}\right)\left(1+\frac{1}{\mathrm{N}(\mathfrak{p})^2}\right).$$

\medskip
$\bullet$ {\bf The Klingen parahoric}
$$U:=\left\{\gamma\in \underline{G}(\hat{\Z}):\,\gamma\equiv\begin{pmatrix}\ast&0&\ast&\ast\\ \ast&\ast&\ast&\ast\\ \ast&0&\ast&\ast\\ 0&0&0&\ast\end{pmatrix}\mod\mathfrak{n}\right\}.$$ 
Consider the module $E$ as above and let $\mathcal{F}_{\underline{G}}(\mathfrak{n})$ be the set of all rank 1 $\left(\mathcal{O}_{F}/\mathfrak{n}\right)$-submodules which are direct factors of $E$. The group $\underline{G}(\hat{\Z})$ acts transitively on  $\mathcal{F}_{\underline{G}}(\mathfrak{n})$ and the stabilizer of the line $L=\langle e_1\rangle$ is the Klingen parahoric $U$. This means that we can identify $\mathcal{F}_{\underline{G}}(\mathfrak{n})$ with the projective space $\mathbf{P}^3(\mathcal{O}_{F}/\mathfrak{n})$.  Thus, its cardinality is given by
$$\#\mathcal{F}_{\underline{G}}(\mathfrak{n})=\mathrm{N}(\mathfrak{n})^3\prod_{\mathfrak{p}\mid\mathfrak{n}}\left(1+\frac{1}{\mathrm{N}(\mathfrak{p})}\right)\left(1+\frac{1}{\mathrm{N}(\mathfrak{p})^2}\right).$$

\medskip
$\bullet$ {\bf The standard Borel parahoric}
$$U:=\left\{\gamma\in \underline{G}(\hat{\Z}):\,\gamma\equiv\begin{pmatrix}\ast&0&\ast&\ast\\ \ast&\ast&\ast&\ast\\ 0&0&\ast& \ast\\ 0&0&0&\ast\end{pmatrix}\mod\mathfrak{n}\right\}.$$ 
Again, we let $E$ be the module as above. We let $\mathcal{F}_{\underline{G}}(\mathfrak{n})$ be the set of all pairs $(L, M)$ such that $L\subset M$ are rank 1 and 2 $\left(\mathcal{O}_{F}/\mathfrak{n}\right)$-submodules of $E$ which are isotropic direct factors. The group $\underline{G}(\hat{\Z})$ acts transitively on  $\mathcal{F}_{\underline{G}}(\mathfrak{n})$ and the stabilizer of the flag $(\langle e_1\rangle, \langle e_1,e_2\rangle)$ is the standard Borel parahoric $U$.

\begin{lem}\label{lem:borel-parahoric}
The cardinality of $\mathcal{F}_{\underline{G}}(\mathfrak{n})$ is given by
$$\#\mathcal{F}_{\underline{G}}(\mathfrak{n})=\mathrm{N}(\mathfrak{n})^4\prod_{\mathfrak{p}\mid\mathfrak{n}}\left(1+\frac{1}{\mathrm{N}(\mathfrak{p})}\right)^2\left(1+\frac{1}{\mathrm{N}(\mathfrak{p})^2}\right).$$
\end{lem}

\begin{prf} Let $M$ be an $(\mathcal{O}_{F}/\mathfrak{n})$-rational point on the Lagrange scheme and $e_1',\,e_2'$ a basis of $M$. A line $L=\langle a_1e_1'+a_2e_2'\rangle$, with $a_1,a_2\in \mathcal{O}_{F}/\mathfrak{n}$, is a direct factor of $M$ if and only if $(a_1,\,a_2)$ determines a point on $\mathbf{P}^1(\mathcal{O}_{F}/\mathfrak{n})$. Combining this with the formula for the cardinality for the Lagrange scheme, we get the lemma by recalling that 
$$\#\mathbf{P}^1(\mathcal{O}_{F}/\mathfrak{n})=\mathrm{N}(\mathfrak{n})\prod_{\mathfrak{p}\mid\mathfrak{n}}\left(1+\frac{1}{\mathrm{N}(\mathfrak{p})}\right).$$
\end{prf}

\begin{rem}\rm It is useful to observe that the proof of Lemma~\ref{lem:borel-parahoric} gives an effective way of listing elements of the flag space of the Borel parahoric. Indeed, by using the description of the Lagrange scheme, we can efficiently list its elements. Then for each such element, we list all the lines that are direct factors by simply running through $\mathbf{P}^1(\mathcal{O}_{F}/\mathfrak{n})$.
\end{rem}

\begin{rem}\rm Lansky and Pollack~\cite{lansky} already make use of the flag space $\mathcal{F}_{\underline{G}}(U)$ in order to find coset representatives for the double coset space $G(\Q)\backslash G(\hat{\Q})/U$. But, they then return to the adelic setting for the computation of the Brandt matrices. One of the advantages of our approach is that, using Proposition~\ref{prop1}, we define the Hecke action on $\mathcal{F}_{\underline{G}}(U)$ directly. 
\end{rem}

\subsection{A sketch of the algorithm}

Here, we give a sketch of the implementation of the algorithm.

\medskip
\begin{alg}[Precomputation]\label{precomputation}\rm The input is a field $F$ as above, a finite set of finite places $\Sigma$ such that $\#\Sigma+[F:\Q]$ is even, and an integral ideal $\mathfrak{n}$ not supported in $\Sigma$.

\begin{enumerate}
\item 
Find a quaternion algebra $D/F$ ramified at precisely the infinite places and the places in $\Sigma$, and compute a maximal order $\CO_D$ of $D$.
\item
Find a set $\Sigma'$ of prime ideals not dividing $\mathfrak{n}$ that generate $\mathrm{Cl}^{+}(F)$ and such that $\Sigma\cap\Sigma'=\emptyset$. 
Choose a prime $\gq\notin\Sigma\cup\Sigma'\cup\{\gp:\gp\mid\gn\}$.
\item 
Compute a complete set $\mathrm{Cl}(\underline{G})$ of representatives $(\mathfrak{a},\alpha)$ for the classes in the genus of $\underline{G}$ which are supported in $(\Sigma',\{\gq\})$ (see Algorithm~\ref{alg:classes-rep}).
\item
For each representative $(\mathfrak{a},\alpha) \in \mathrm{Cl}(\underline{G})$, compute the group $\overline{\Gamma}_{\mathfrak{a},\alpha}$.
\end{enumerate}
\end{alg}

\medskip
\begin{alg}[Main algorithm]\label{main-alg}\rm The input consists of $F$, the level $\mathfrak{n}$ (with the level structure specified: Borel, Klingen or Siegel) and the precomputed data, together with a weight representation $(\rho,V)$. The output consists of a set of primes $\mathfrak{p}$ (ordered by norm), the Hecke matrices $T_1(\mathfrak{p})$ and $T_2(\mathfrak{p})$ and the Hecke constituents.

\begin{enumerate}
\item Compute splitting isomorphisms $\GU_2(\CO_{D_\gq})\simeq\GSp_4(\mathcal{O}_{F_\gq})$, for each prime $\gq\mid\mathfrak{n}$.
\item For each $(\mathfrak{a},\alpha)\in \mathrm{Cl}(\underline{G})$, compute $M_{\rho}^{\Gamma_{\mathfrak{a},\alpha}}(\mathfrak{n})$ as a module of coinvariants
$$M_{\rho}^{\Gamma_{\mathfrak{a},\alpha}}(\mathfrak{n})=E[\mathcal{F}_{\underline{G}}(\mathfrak{n})]\otimes V/\langle x-\gamma x,\,\gamma\in\Gamma_{\mathfrak{a},\alpha}\rangle.$$
\item Combine the results of step (2), forming the direct sum $$M_{\rho}(\mathfrak{n})=\bigoplus_{(\mathfrak{a},\alpha)\in \mathrm{Cl}(\underline{G})}M_{\rho}^{\Gamma_{\mathfrak{a},\alpha}}(\mathfrak{n}).$$
\item Starting with the primes of smallest norm in $\CO_F$, compute the Hecke operators $T_1(\mathfrak{p})$ and $T_2(\mathfrak{p})$ for all the primes $\gp$, in increasing order, until you find a common basis of eigenvectors of $M_{\rho}(\mathfrak{n})$ for the $T_1(\mathfrak{p})$ and $T_2(\mathfrak{p})$ that completely diagonalize $M_{\rho}(\mathfrak{n})$. 
\end{enumerate}
\end{alg}

\begin{rem}\rm In practice, one can improve the efficiency of Step (4) of Algorithm~\ref{main-alg}, by diagonalizing $M_{\rho}(\gn)$ only using the Hecke operators $T_1(\gp)$ and then computing the $T_2(\gp)$ when needed.
\end{rem}

 \section{Examples of Hilbert-Siegel eigenforms on $\Q(\sqrt{2})$}\label{sec:example}
 In this section, we give some examples of Hilbert-Siegel eigensystems for $F=\Q(\sqrt{2})$. (We obtained them with a preliminary version of our algorithm, which we implemented in {\sf Magma}~\cite{magma}.) 
 
\subsection{Hilbert-Siegel eigensystems of level 1}\label{subsec:level-1}
Let $D$ be the Hamilton quaternion algebra over $F$, i.e., the quaternion algebra over $F$ determined by the relations $i^2=-1$, $j^2=-1$ and $k=ij=-ji$. Then $\mathrm{disc}(D)=(1)$ and we choose the maximal order $\CO_D=\Z[\sqrt{2}][e_1,\,e_2,\,e_3,\,e_4]$, where
\begin{eqnarray*}
e_1=\frac{1+i}{\sqrt{2}},\,\,e_2=\frac{1+j}{\sqrt{2}},\,e_3=e_1e_2,\,e_4=e_2e_1.
\end{eqnarray*}
We recall that the narrow class number of $F$ and the class number of the quaternion algebra $D$ are both equal to 1. Thus every representative in the principal genus of $G$ is of the form $((1),\alpha)$, where $\alpha\bar{\alpha}^t=\gamma\in\GL_2(\CO_D)$. We choose $\Sigma_U=(\emptyset,\{(\sqrt{2})\})$. We determine that the matrices
$$\gamma_1=\begin{pmatrix}1&0\\ 0&1 \end{pmatrix},\,\gamma_2=\begin{pmatrix}2& \frac{(2+\sqrt{2}) + (2-\sqrt{2})i}{2}\\ \frac{2+\sqrt{2} + (- 2+\sqrt{2})i}{2}& 2\end{pmatrix} $$
correspond to distinct equivalence classes by computing their stabilizers for which we obtain $\mathrm{Card}(\Gamma_1)=4608$ and $\mathrm{Card}(\Gamma_2)=3840$, where we simply let $\Gamma_i=\overline{\Gamma}_{(1),\alpha_i}$. By using the mass formula in W. T. Gan, J. Hanke and J.-K. Yu \cite[Proposition 9.3]{gan1}, we verify that
\begin{eqnarray*}
\mathrm{Mass}(\underline{G})&=&\frac{1}{2^4}\zeta_{\Q(\sqrt{2})}(-1)\zeta_{\Q(\sqrt{2})}(-3)=\frac{1}{2^4}\frac{1}{12}\frac{11}{120}=\frac{1}{4608}+\frac{1}{3840}\\
&=&\frac{1}{\mathrm{Card}(\Gamma_1)}+\frac{1}{\mathrm{Card}(\Gamma_2)}.
\end{eqnarray*}  Therefore, the class number of the unique genus of $G$ is $h=2$. 

This means that the dimension of the space of holomorphic Hilbert-Siegel modular forms of level $1$ and parallel weight $3$ is $2$. From this, we determine that there is only one Hilbert-Siegel cuspidal eigenform $f$ of level $1$ and parallel weight 3. The Brandt matrices $\mathcal{B}_1(\mathfrak{p})$ and $\mathcal{B}_2(\mathfrak{p})$ of the Hecke operators acting on $M_{3}(1)$ for the first four prime ideals $\mathfrak{p}$ are given below.
\begin{eqnarray*}
\begin{array}{cc|cc}
\mathrm{N}(\mathfrak{p})&\mathfrak{p}&\mathcal{B}_1(\mathfrak{p})&\mathcal{B}_2(\mathfrak{p})\\\hline
&&&\\
2&(\sqrt{2})&\begin{pmatrix} 9& 6\\ 5 &10\end{pmatrix}&\begin{pmatrix}12 &18\\ 15 &15\end{pmatrix}\\
&&&\\
7&(3+\sqrt{2})&\begin{pmatrix}208 &192\\ 160& 240\end{pmatrix}&\begin{pmatrix}1264 &1536\\ 1280 &1520\end{pmatrix}\\
&&&\\
7&(3-\sqrt{2})&\begin{pmatrix}208 &192\\ 160& 240\end{pmatrix} &\begin{pmatrix}1264 &1536\\ 1280 &1520\end{pmatrix}\\
&&&\\
9&(3)&\begin{pmatrix} 436 &384\\ 320 &500\end{pmatrix}&\begin{pmatrix} 3540 &3840\\ 3200& 4180\end{pmatrix}
\end{array}
\end{eqnarray*}
The first few Hecke eigenvalues of the eigenform $f$ are listed in Table~\ref{table: table1}.  The eigenvalues of $f$ (or in fact the Brandt matrices) at the primes above $7$ suggest that it is a base change from $\Q$, which we then confirm as follows. For a classical or Hilbert newform $h$ of weight $4$, let $SK(h)$ be its Saito-Kurokawa lift.  Letting $a_\gp(h)$ be the Hecke eigenvalue of $h$ at $\gp$, the Hecke eigenvalues $\lambda_1(\gp)$ and $\lambda_2(\gp)$ of $SK(h)$ corresponding to the Hecke operators $T_1(\gp)$ and $T_2(\gp)$ are then given by
\[
\begin{aligned}
\lambda_1(\mathfrak{p})&=a_\gp(h)\ \mathrm{N}(\mathfrak{p})^{\frac{4-k}{2}}+\mathrm{N}(\mathfrak{p})^2+\mathrm{N}(\mathfrak{p})\\
\lambda_2(\mathfrak{p})&=a_\gp(h)\ \mathrm{N}(\mathfrak{p})^{\frac{4-k}{2}}(\mathrm{N}(\mathfrak{p})+1)+\mathrm{N}(\mathfrak{p})^2-1.
\end{aligned}
\]
\begin{prop}\label{prop:base-change} The form $f$ is a base change of the form $SK(h)$, where $h$ is a newform which belongs to the unique conjugacy class of $S_4(8,(\frac{\cdot}{F}))^{\mbox{\tiny\rm new}}$, the newspace of classical forms of level 8 and character $(\frac{\cdot}{F})$.
\end{prop}

\begin{prf}
We compute the space $S_4(1)$ of Hilbert cusp forms of level 1 and pa\-ral\-lel weight $4$ over $F$ using the Hilbert Modular Forms Package in {\sf Magma}~\cite{magma}, and obtain that it has dimension 1. We then identify the eigenform $f$ as the Saito-Kurokawa lift of the newform $g$ in $S_4(1)$ by direct calculations. By Saito~\cite[Theorem 4.5]{saito1}, the form $g$ itself is a base change of a newform $h$ in $S_4(8,(\frac{\cdot}{F}))^{\mbox{\tiny\rm new}}$. The lift $SK(h)$ is a classical Siegel eigenform of level 8 and character $(\frac{\cdot}{F})$. By the funtoriality of Saito-Kurokawa lifts~\cite{cogdell-shapiro}, we get
$$f=SK(g)=SK(BC_{F/\Q}(h))=BC_{F/\Q}(SK(h)).$$
\end{prf}

\begin{rem}\rm Proposition~\ref{prop:base-change}, combined with Sorensen~\cite[Theorem B]{sorensen1}, implies that there is no holomorphic Hilbert-Siegel eigenform of level 1 and parallel weight 3 over $\Q(\sqrt{2})$ which is stable.
\end{rem}

By lifting of Galois representations attached to $h$, there is a family of $\ell$-adic Galois representations $(\rho_{SK(h),\ell})$  attached to $SK(h)$:
$$\rho_{SK(h),\ell}:\,\Gal(\overline{\Q}/\Q)\to\mathrm{GSp}_4(\overline{\Q}_\ell).$$
Similarly, there is a family $(\rho_{f,\ell})$ of $\ell$-adic representations of $\Gal(\overline{\Q}/F)$ attached to $f$, which is everywhere unramified since $f$ has level $1$. In fact, by Proposition~\ref{prop:base-change}, we simply have $\rho_{f,\ell}=\rho_{SK(h),\ell}|_{\Gal(\overline{\Q}/F)}$.

It would be interesting to know whether there exists a (Calabi-Yau) threefold $X/\Q(\sqrt{2})$, with Hodge numbers $h^{3,0}=h^{2,1}=h^{1,2}=h^{0,3}=1$, such that for every prime $\ell$ the Galois representation $$\rho_{X,\ell}:\,\Gal(\overline{\Q}/\Q(\sqrt{2}))\to\GL(H_{\mbox{\tiny\'et}}^3(X,\,\overline{\Q}_\ell))$$ is isomorphic to $\rho_{f,\ell}$. The existence of such a threefold would provide an analogue in higher rank of the Tate (or Shimura) elliptic curve $E/\Q(\sqrt{29}): y^2+xy+\varepsilon^2 y=x^3$ with discriminant $\Delta_E=-\varepsilon^{10}$, where $\varepsilon=\frac{5+\sqrt{29}}{2}$ is the fundamental unit, in the sense that $\Q(\sqrt{2})$ would be the totally real field of smallest discriminant for which there is a higher rank (automorphic) motive over $\Q$ whose base change to $\Q(\sqrt{2})$ becomes everywhere unramified. (Indeed, the results in~\cite{clifton1} show that no such threefold exists over $\Q(\sqrt{5})$.)

\subsection{Hilbert-Siegel eigensystems of prime level}\label{subsec:prime-levels}

For all prime levels $\mathfrak{r}$, with odd norm $\mathrm{N}(\mathfrak{r})\le 31$ (and up to Galois conjugation), we compute the space $S_{\rho}(\mathfrak{r})$ using the same precomputed data as in~\ref{subsec:level-1}, for the trivial representation $\rho$ and the Siegel parahoric. By the Jacquet-Langlands correspondence for $\GSp_4$ for fields of even degree~\cite[Theorem B]{sorensen1}, we obtain Hecke eigensystems of Hilbert-Siegel cuspidal forms of parallel weight 3. (Note that \cite[Theorem B]{sorensen1} is only true for eigensystems which come from stable automorphic forms; but for forms that are lifts, one can use the usual Jacquet-Langlands correspondence for $\GL_2$.) In Tables~\ref{table: table1} and~\ref{table: table2}, we list all those eigensystems that are defined over $\Q$ or a quadratic field. Here are the conventions we use in the tables.
\begin{enumerate}
\item For a quadratic field $K$ of discriminant $d$, we let $\omega_d$ be a generator of the ring of integers $\mathcal{O}_K$ of $K$.
\item The first row contains the level $\mathfrak{r}$ given in the format $(\mathrm{N}(\mathfrak{r}), \alpha)$, where $\alpha\in\mathcal{O}_F$ is a generator of $\mathfrak{r}$, and the dimensions of the relevant spaces. 
\item The second and third rows list the computed Hecke operators.
 \item The Hecke eigenvalues systems are given by row.
\end{enumerate} 
 
\begin{landscape}
\begin{table}
\begin{eqnarray*}
\begin{array}{|crrrrrrrr|}\hline
\multicolumn{9}{|l|}{\mathfrak{r}=(1,\,1),\quad\dim M_3(\mathfrak{r})=2,\quad\dim S_3(\mathfrak{r})=1}\\\hline
\mathrm{N}(\mathfrak{p})&\multicolumn{2}{c}{2}&\multicolumn{2}{c}{7}&\multicolumn{2}{c}{7}&\multicolumn{2}{c|}{9}\\\hline
&T_1(2+\omega_{8})&T_2(2+\omega_{8})&T_1(3+\omega_{8})&T_2(3+\omega_{8})&T_1(3-\omega_{8})&T_2(3-\omega_{8})&T_1(3)&T_2(3)\\\hline\hline
f_{1}&4&-3&48&-16&48&-16&116&340\\
\hline
\multicolumn{9}{c}{}\\\hline
\multicolumn{9}{|l|}{\mathfrak{r}=(7,\,3+\omega_{8}),\quad\dim M_3(\mathfrak{r})=6,\quad\dim S_3(\mathfrak{r})=5}\\\hline
\mathrm{N}(\mathfrak{p})&\multicolumn{2}{c}{2}&\multicolumn{2}{c}{7}&\multicolumn{2}{c}{7}&\multicolumn{2}{c|}{9}\\\hline
&T_1(2+\omega_{8})&T_2(2+\omega_{8})&T_1(3+\omega_{8})&T_2(3+\omega_{8})&T_1(3-\omega_{8})&T_2(3-\omega_{8})&T_1(3)&T_2(3)\\\hline\hline
f_{1}&10&15&-7&0&60&80&80&-20\\
f_{2}&-4&1&7&0&32&80&-60&148\\
\hline
\multicolumn{9}{c}{}\\\hline
\multicolumn{9}{|l|}{\mathfrak{r}=(9,\,3),\quad\dim M_3(\mathfrak{r})=12,\quad\dim S_3(\mathfrak{r})=11}\\\hline
\mathrm{N}(\mathfrak{p})&\multicolumn{2}{c}{2}&\multicolumn{2}{c}{7}&\multicolumn{2}{c}{7}&\multicolumn{2}{c|}{9}\\\hline
&T_1(2+\omega_{8})&T_2(2+\omega_{8})&T_1(3+\omega_{8})&T_2(3+\omega_{8})&T_1(3-\omega_{8})&T_2(3-\omega_{8})&T_1(3)&T_2(3)\\\hline\hline
f_{1}&-6&7&-22&64&-22&64&-9&0\\
f_{2}&4&5&-16&16&-16&16&9&0\\
\hline
\multicolumn{9}{c}{}\\\hline
\multicolumn{9}{|l|}{\mathfrak{r}=(17,\,5+2\omega_{8}),\quad\dim M_3(\mathfrak{r})=23,\quad\dim S_3(\mathfrak{r})=22}\\\hline
\mathrm{N}(\mathfrak{p})&\multicolumn{2}{c}{2}&\multicolumn{2}{c}{7}&\multicolumn{2}{c}{7}&\multicolumn{2}{c|}{9}\\\hline
&T_1(2+\omega_{8})&T_2(2+\omega_{8})&T_1(3+\omega_{8})&T_2(3+\omega_{8})&T_1(3-\omega_{8})&T_2(3-\omega_{8})&T_1(3)&T_2(3)\\\hline\hline
f_{1}&-2&3&6&32&-36&80&8&28\\
f_{2}&4-\omega_{40}&-3-3\omega_{40}&38-4\omega_{40}&-96-32\omega_{40}&58+6\omega_{40}&64+48\omega_{40}&66+4\omega_{40}&-160+40\omega_{40}\\\hline
\end{array}\end{eqnarray*}
\caption{\bf Hilbert-Siegel eigenforms of parallel weight 3 over $\Q(\sqrt{2})$}
\label{table: table1}
\end{table}
\end{landscape}

\begin{landscape}
\begin{table}
\begin{eqnarray*}
\begin{array}{|crrrrrrrr|}\hline
%\multicolumn{9}{c}{}\\\hline
\multicolumn{9}{|l|}{\mathfrak{r}=(23,\,5+\omega_{8}),\quad\dim M_3(\mathfrak{r})=32,\quad\dim S_3(\mathfrak{r})=31}\\\hline
\mathrm{N}(\mathfrak{p})&\multicolumn{2}{c}{2}&\multicolumn{2}{c}{7}&\multicolumn{2}{c}{7}&\multicolumn{2}{c|}{9}\\\hline
&T_1(2+\omega_{8})&T_2(2+\omega_{8})&T_1(3+\omega_{8})&T_2(3+\omega_{8})&T_1(3-\omega_{8})&T_2(3-\omega_{8})&T_1(3)&T_2(3)\\\hline\hline
f_{1}&-8&10&-46&119&-54&145&-54&153\\
f_{2}&10&15&72&176&56&48&124&420\\
\hline
\multicolumn{9}{c}{}\\\hline
\multicolumn{9}{|l|}{\mathfrak{r}=(25,\,5),\quad\dim M_3(\mathfrak{r})=48,\quad\dim S_3(\mathfrak{r})=47}\\\hline
\mathrm{N}(\mathfrak{p})&\multicolumn{2}{c}{2}&\multicolumn{2}{c}{7}&\multicolumn{2}{c}{7}&\multicolumn{2}{c|}{9}\\\hline
&T_1(2+\omega_{8})&T_2(2+\omega_{8})&T_1(3+\omega_{8})&T_2(3+\omega_{8})&T_1(3-\omega_{8})&T_2(3-\omega_{8})&T_1(3)&T_2(3)\\\hline\hline
f_{1}&2&-9&62&96&62&96&40&-420\\
f_{2}&-4-2\omega_{12}&5+2\omega_{12}&-15+7\omega_{12}&56-8\omega_{12}&-15+7\omega_{12}&56-8\omega_{12}&44-18\omega_{12}&132-52\omega_{12}\\
f_{3}&8&9&44-10\omega_{24}&-48-80\omega_{24}&44+10\omega_{24}&-48+80\omega_{24}&76&-60\\
\hline
\multicolumn{9}{c}{}\\\hline
\multicolumn{9}{|l|}{\mathfrak{r}=(31,\,7+3\omega_{8}),\quad\dim M_3(\mathfrak{r})=65,\quad\dim S_3(\mathfrak{r})=64}\\\hline
\mathrm{N}(\mathfrak{p})&\multicolumn{2}{c}{2}&\multicolumn{2}{c}{7}&\multicolumn{2}{c}{7}&\multicolumn{2}{c|}{9}\\\hline
&T_1(2+\omega_{8})&T_2(2+\omega_{8})&T_1(3+\omega_{8})&T_2(3+\omega_{8})&T_1(3-\omega_{8})&T_2(3-\omega_{8})&T_1(3)&T_2(3)\\\hline\hline
f_{1}&-4&4&-60&174&-20&-2&-60&130\\
f_{2}&-4&5&20&16&-8&48&8&28\\
f_{3}&-6&7&44-2\omega_{204}&112-8\omega_{204}&-4\omega_{204}&48&-32-2\omega_{204}&108+4\omega_{204}\\
f_{4}&2&-9&72-2\omega_{204}&176-16\omega_{204}&56-4\omega_{204}&48-32\omega_{204}&76-2\omega_{204}&-60-20\omega_{204}\\
\hline\end{array}\end{eqnarray*}
\caption{\bf Hilbert-Siegel eigenforms of parallel weight 3 over $\Q(\sqrt{2})$ (cont'd)}
\label{table: table2}
\end{table}
\end{landscape}


\begin{thebibliography}{99}
\bibitem{bellaiche1} J.~Bella\"{\i}che and P.~Graftieaux, Augmentation du niveau pour ${\rm U}(3)$. {\it Amer. J. Math}. {\bf 128} (2006), no. 2, 271--309.

\bibitem{magma} 
W.~Bosma, J.~Cannon, and C.~Playoust, \emph{The Magma algebra system. I. The user language}, J.~Symbolic Comput. \textbf{24} (1997), vol.~3--4, 235--265.

\bibitem{cogdell-shapiro} J. W. Cogdell and I. I.~Piatetski-Shapiro, Base change for the Saito-Kurokawa representations of ${\rm PGSp}(4)$. {\it J. Number Theory} {\bf 30} (1988), no. 3, 298--320.

\bibitem{coulangeon} R. Coulangeon, Tensor products of Hermitian lattices. {\it Acta Arith}. {\bf 92} (2000), no. 2, 115--130.

\bibitem{clifton1} C. Cunningham and L. Demb\'el\'e, Computing genus 2 Hilbert-Siegel modular forms  over $\Q(\sqrt{5})$ via the Jacquet-Langlands Correspondence (14 pages). To appear in {\it Experimental Mathematics}.

\bibitem{dembele2} L. Demb\'el\'e, Quaternionic $M$-symbols, Brandt matrices and Hilbert modular forms. {\it Math. Comp.} Vol. {\bf 76}, no 258 (2007), 1039-1057.

\bibitem{dd1} L. Demb\'el\'e and S. Donnelly, Computing Hilbert modular forms over fields with nontrivial class group. Algorithmic number theory, 371--386, {\it Lecture Notes in Comput. Sci}., {\bf 5011}, Springer, Berlin, 2008.

\bibitem{eichler1} M.~Eichler, On theta functions of real algebraic number fields. {\it Acta Arith}. {\bf 33} (1977), no. 3, 269--292.

\bibitem{gan1} W. T. Gan, J. P.~Hanke, and J.-K.~Yu, On an exact mass formula of Shimura. {\it Duke Math. J.} {\bf 107} (2001), no. 1, 103--133.

\bibitem{ghitza1}
A.~Ghitza, Hecke eigenvalues of Siegel modular forms (mod $p$) and of algebraic modular forms. {\it J. Number Theory} \textbf{106} (2004), no. 2, 345--384.

\bibitem{gross1} B. H.~Gross, Algebraic modular forms. {\it Israel J. Math.} {\bf 113} (1999), 61--93.

\bibitem{gross2} B. H.~Gross, Groups over $\Z$. {\it Invent. Math.} {\bf 124} (1996), no. 1-3, 263--279.

\bibitem{hashimoto1} K. Hashimoto and T. Ibukiyama, On class numbers of positive definite binary quaternion Hermitian forms. II. {\it J. Fac. Sci. Univ. Tokyo Sect. IA Math.} {\bf 28} (1981), no. 3, 695--699. 

\bibitem{HerzigTilouine} F.~Herzig and J.~Tilouine, Conjecture de type de Serre et formes compagnons pour $\GSp_4$ (preprint). Available at:
{\sf http://math.northwestern.edu/~herzig/ht08-7c.pdf}.

\bibitem{ibukiyama}  T.~Ibukiyama, On symplectic Euler factors of genus two. {\it J. Fac. Sci. Univ. Tokyo Sect. IA Math.} {\bf 30} (1984), no. 3, 587--614.

\bibitem{ihara} Y. Ihara, On certain arithmetical Dirichlet series. {\it J. Math. Soc. Japan} Vol. {\bf 16} no 3 (1964), 214--225.

\bibitem{lansky} J.~Lansky, and D.~Pollack, Hecke algebras and automorphic forms. {\it Compositio Math.} {\bf 130} (2002), no. 1, 21--48.

\bibitem{ramakrishnan1} D.~Ramakrishnan and F.~Shahidi, Siegel modular forms of genus 2 attached to elliptic curves. {\it Math. Res. Lett}. {\bf 14} (2007), no. 2, 315--332.

\bibitem{reiner} I. Reiner, {\it Maximal orders}. London Mathematical Society Monographs, No. 5. Academic Press, London-New York, 1975. xii+395 pp.

\bibitem{saito1} H. Saito, Automorphic forms and algebraic extensions of number fields. II. {\it J. Math. Kyoto Univ.} {\bf 19} (1979), no. 1, 105--123.

\bibitem{shimura1} G. Shimura, Arithmetic of unitary groups. {\it Ann. of Math.} (2) {\bf 79} (1964) 369--409.

\bibitem{shimura2} G. Shimura, Arithmetic of alternating forms and quaternion hermitian forms. {\it J. Math. Soc. Japan} {\bf 15} (1963), 33--65.

\bibitem{sorensen1} C.~M.~Sorensen, Potential level-lowering for $\GSp(4)$. {\it J. Inst. Math. Jussieu} {\bf 8}, 595-622. 

%\bibitem{weissauer1} R.~Weissauer : {\it Four dimensional Galois representations}, in Formes Automorphes (II), le cas du groupe $\GSp(4)$, Ast\'erisque 302, pp. 67Ð150, 2005. 

\end{thebibliography}
\end{document}